\documentclass[11pt]{amsart}

\usepackage{amsmath,amsthm, amscd, amssymb, amsfonts}
\usepackage[all]{xy}
\usepackage[dvips, dvipsnames, usenames]{color}

\newtheorem{theorem}{Theorem}[section]

\newtheorem{lem}[theorem]{Lemma}

\newtheorem{pro}[theorem]{Proposition}

\theoremstyle{definition}

\theoremstyle{remark}

\newcommand{\ydh}{{}^{H}_{H}\mathcal{YD}}
\newcommand{\ydg}{{}^{\Gamma}_{\Gamma}\mathcal{YD}}

\newcommand{\VGamma}{\widehat{\Gamma}}
\newcommand{\e}{\mathbf e}

\def\zt{\Z^{\theta}}

\newcommand{\otvz}{{1\le i, j \le \theta}}
\def\G{\mathbb{G}}
\newcommand\id{\operatorname{id}}

\newcommand\ord{\operatorname{ord}}
\newcommand\Hom{\operatorname{Hom}}

\newcommand\Aut{\operatorname{Aut}}

\newcommand\gr{\operatorname{gr}}

\newcommand\ad{\operatorname{ad}}

\def\gyd{{}^{\k \Gamma}_{\k \Gamma}\mathcal{YD}}

\def\k{\Bbbk}

\def\ot{\otimes}

\def\s{\mathbb{S}}

\def\Z{\mathbb{Z}}
\def\N{\mathbb{N}}
\def\X{\mathbb{X}}
\def\B{\mathfrak{B}}

\def\eps{\epsilon}

\def\mD{\mathcal{D}}

\def\mP{\mathcal{P}}

\def\bS{\mathfrak{S}}
\def\bB{\mathfrak{B}}

\def\pf{\begin{proof}}
\def\epf{\end{proof}}

\begin{document}


\title[Generation in degree one]
{Pointed Hopf algebras with standard braiding are generated in degree one}
\author[Angiono; Garc\'ia Iglesias]
{Angiono, Iv\'an; Garc\'ia Iglesias, Agust\'in}

\address{FaMAF-CIEM (CONICET), Universidad Nacional de C\'ordoba,
Medina A\-llen\-de s/n, Ciudad Universitaria, 5000 C\' ordoba, Rep\'
ublica Argentina.} \email{(angiono|aigarcia)@mate.uncor.edu}

\thanks{The work was partially supported by CONICET,
FONCyT-ANPCyT, Secyt (UNC), Mincyt (C\'ordoba)}

\begin{abstract}
We show that any finite-dimensional pointed Hopf algebra over an abelian group $\Gamma$ such that its infinitesimal braiding is of standard type is generated by group-like and skew-primitive elements. This fact agrees with the long-standing conjecture by Andruskiewitsch and Schneider. We also show that the quantum Serre relations hold in any coradically graded pointed Hopf algebra over $\Gamma$ of finite dimension and determine how these relations are lifted in the standard case.
\end{abstract}

\maketitle

\section*{Introduction}

The classification of finite-dimensional pointed Hopf algebras is currently an active area of research. This class includes group algebras $\k \Gamma$, $\Gamma$ a group, and Frobenius-Luzstig kernels $\mathfrak{u}_q(\mathfrak{g})$ \cite{L}, associated to a semisimple finite-dimensional Lie algebra $\mathfrak{g}$ and a root of unity $q$.

The main result concerning the classification of finite-dimensional pointed Hopf algebras with group of group-likes $\Gamma$ over an algebraically closed field of characteristic zero has been done by Andruskiewitsch and Schneider in \cite{AS4}, for the case $\Gamma$ abelian and $|\Gamma|$ not divisible by $2,3,5,7$. This was achieved using the so-called Lifting method, introduced by the authors in previous works, see e.g. \cite{AS1}, \cite{AS3}.

One of the main steps in the Lifting method is to determine when a given braiding yields a finite-dimensional Nichols algebra. This was solved in the abelian case by Heckenberger \cite{H2}. Another key step of this method is to prove that all pointed Hopf algebras over $\Gamma$ are generated by group-like and skew-primitive elements, or, equivalently, that the associated graded Hopf algebra with respect to the coradical filtration is generated in degrees 0 and 1. This problem has been solved for a finite-dimensional pointed Hopf algebra $H$ over a group $\Gamma$ in the following cases:
\begin{itemize}
\item when $H$ is co-triangular \cite{AEG},
\item when $\Gamma$ is abelian and $|\Gamma|$ is not divisible by $2,3,5,7$ \cite{AS4},
\item when $\Gamma$ is (isomorphic to) $\s_n$, $n=3,4,5$ \cite{AG, GG},
\item when the braiding arises from some particular affine racks \cite{AG}.
\end{itemize}
It has been conjectured in \cite[Conjecture 1.4]{AS2} that this holds for any $H$ as above. Our main result, Theorem \ref{thm:generationdegree1} is a positive answer to this conjecture in the case in which the braiding is of standard type.

The paper is organized as follows. In Section 1 we recall the basic facts about braided vector spaces and Nichols algebras, including the description of standard braidings. Also we explicitly recall the steps of the Lifting method.

In Section 2 we prove Theorem \ref{thm:generationdegree1}. To do this, we first prove in Proposition \ref{pro:qSerre} that quantum Serre relations hold in any finite-dimensional graded braided Hopf algebra in the category of Yetter-Drinfeld modules over an abelian group $\Gamma$, not necessarily of standard type. To complete the proof we strongly use the presentation of Nichols algebras of standard type by generators and relations given in \cite{A}, which we recall in Theorem \ref{thm:presentacionNichols}. Along this part, we repeatedly use the classification of finite arithmetic root systems in \cite{H2} to check that some diagrams that we can associate to the relations are not of finite type.

In Section 3 we show how the quantum Serre relations are lifted to a general pointed Hopf algebra, exploiting the arguments in the proof of Proposition \ref{pro:qSerre}.

\section{Preliminaries}

\subsection{Conventions}
We work over an algebraically closed field $\k$ of characteristic $0$.
For each $N > 0$, $\G_N$ denotes the group of primitive $N$-th roots
of 1 in $\k$. Given $n \in \N$ and $q \in \k$, $q \notin \cup_{0 \leq j \leq n}
\G_j$, $q\neq0$, we denote
$$ \binom{n}{j}_q = \frac{(n)_q!}{(k)_q! (n-k)_q!}, \quad \mbox{where }(n)_q!= \prod_{j=1}^n (k)_q \quad \mbox{and } (k)_q= \sum_{j=0}^{k-1} q^j. $$

If $\Gamma$ is an abelian group, we denote $\widehat{\Gamma}= \Hom(\Gamma, \k^{\times})$.

\subsection{Nichols algebras and the Lifting method}

\subsubsection{Braided vector spaces}
A braided vector space is a pair $(V,c)$, where $V$ is a vector
space and $c\in \Aut (V\ot V)$ is a solution of the braid
equation: $$(c\otimes \id) (\id\otimes c) (c\otimes \id) =
(\id\otimes c) (c\otimes \id) (\id\otimes c).$$ We extend the
braiding to $c:T(V)\ot T(V) \to T(V)\ot T(V)$ in the usual way. If
$x,y\in T(V)$, then the braided commutator is
\begin{equation*}\label{eqn:braidedcommutator}
[x,y]_c := \text{multiplication } \circ \left( \id - c \right)
\left( x \ot y \right).
\end{equation*}

\subsubsection{Braided vector spaces of diagonal type and hyperwords}\label{subsec:braideddiagonal}
A braided vector space $(V,c)$ is of \emph{diagonal type} with
respect to a basis $\{x_i\}_{i\in I}$ if there exist
$q_{ij}\in \k^{\times}$ such that $c(x_i \ot x_j)= q _{ij} x_j
\ot x_i$, $i,j\in I$.
\medbreak

Fix a basis $\{x_1,\dots ,x_{\theta}\}$ of $V$ (we assume $\dim V<\infty$). Denote by $\X$ the set of words in letters $x_1, \ldots , x_\theta$ and order them using the lexicographic order. We identify canonically $\X$ with a basis of $T(V)$.

We say that $u\in \X$ is a \emph{Lyndon word} if $u$ is smaller than any of its proper ends. That is, $u=x_i$ for some $i$ or for each $v,w \in \X\setminus\{1\}$ such that $u=vw$, we have $u<w$.

Denote by $L$ the set of Lyndon words. It follows that $u\in L$ if and only if there exist $v<w \in L$ such that $u=vw$ or $u=x_i$, $1\leq i\leq \theta$. Given $u\in L\setminus\{x_1,\ldots,x_\theta\}$ the \emph{Shirshov decomposition} $u=vw$ with $v,w \in L$ is the one such that $w$ is the smallest end between all possible decompositions. See \cite{Kh} and references therein.

For each $u \in L$ we consider an element $[u]_c \in T(V)$, called the \emph{hyperletter} \cite{Kh} corresponding to $u$, defined inductively by
$$[u]_c=\left\{ \begin{array}{ll}
            u, & \hbox{if }u\in \X;  \\
            \left[ [v]_c, [w]_c \right]_c, & \hbox{if }u=vw \hbox{ is the Shirshov decomposition of }u.
          \end{array}  \right. $$

\subsubsection{Braided vector spaces of standard type}\label{subsec:bvs}
Fix $\theta\in\N$, $I=\{1,\ldots,\theta\}$ and $(q_{ij})_\otvz$ as above.

Let $E=\{\e_1,\dots,\e_\theta\}$ be the canonical basis of
$\Z^\theta$ and $\chi: \zt\times \zt \to \k^{\times}$ the bilinear form determined by $\chi(\e_i, \e_j)=q_{ij}$,$\otvz$.  If $F=\{f_1\dots,f_\theta\}$ is another basis of $\Z^\theta$, then
we set $q^F=(q^F_{ij})_{1\leq i,j\leq \theta}$, by
$q^F_{ij}=\chi(f_i,f_j)$, $1\leq i,j\leq \theta$. Thus $q=q^E$. For
$i\neq j\in\{1,\dots,\theta\}$, consider the set
\begin{equation*}
  M^F_{ij}=\{m\in\N_0 \vert
  (m+1)_{q^F_{ii}}({q^F_{ii}}^m q^F_{ij}q^F_{ji}-1)=0\}.
\end{equation*}
If this set is not empty, let $m^F_{ij}$ denote its minimal element.
Also let $m^F_{ii}=2, \forall\,i$. Let $s^F_i$ be the
pseudo-reflection in $\Z^\theta$ given by
\begin{equation*}
  s^F_i(f_j)=f_j+m^F_{ij}f_i, \ j=1,\dots,\theta.
\end{equation*}
Let $\Omega$ be the set of all ordered bases of $\Z^\theta$ and let
$\mP(\chi)\subseteq\Omega$ be the set of points of the Weyl grupoid
$W(\chi)$ of the bilinear form $\chi$, see \cite[Def. 3.2]{A}.

Clearly, $s^F_i(F)$ is again a basis of $\Z^\theta$. The form $\chi$
is called \emph{standard} if for every $F\in\mP(\chi)$, the integers
$m^F_{ij}$ are defined, for all $1\leq i,j\leq \theta$, and
$m(s^F_k(F))_{ij}=m_{ij}$ for every $i,j,k$. The corresponding braided vector space is said to be of \emph{standard type} \cite{AA}.

We are interested in standard braiding whose associated Nichols algebra is finite-dimensional. In such case, the corresponding Cartan matrix $C=(a_{ij}=-m_{ij})_{i,j \in \{1, \ldots, \theta \} }$ is finite, see \cite[Thm. 4.1]{A}. This family includes properly the braidings of Cartan type considered in \cite{AS2}.

Standard braidings with finite-dimensional Nichols algebras are classified in \cite{A}. In the same paper, the dimension, a presentation by generators and relations and a PBW basis is given for each Nichols algebras with standard braiding. We recall next this result, which will be fundamental for our work.

\begin{theorem}{\cite[Thms. 5.14, 5.19, 5.22, 5.25]{A}}
\label{thm:presentacionNichols}
Let $V$ be a braided vector space of standard type, of dimension $\theta$, $C=(a_{ij}=-m_{ij})_{i,j \in \{1, \ldots, \theta \} }$ the associated finite Cartan matrix, and $\Delta_+$ the corresponding root system.

The Nichols algebra ${\mathfrak B}(V)$ is presented by generators
$x_{i}$, $1\le i \le \theta$, and the following relations
\begin{align}
\label{alturaPBW} x_{\alpha}^{N_{\alpha}} &= 0, \quad \alpha \in
\Delta_+;
\\ ad_{c}(x_{k})^{1+m_{kj}}(x_{j})  &= 0, \quad k\neq j, \ q_{kk}^{m_{kj}+1} \neq 1;
\label{serrebis}
\end{align}
if there exist distinct $j, k, l$ satisfying $\ m_{kj}=m_{kl}=1$,
$q_{kk}=-1$, then
\begin{equation} \left[ (\ad x_k)x_j, (\ad
x_k)x_l \right]_c = 0;\label{Arelation}
\end{equation}
if there exist $k \neq j$ satisfying $\ m_{kj}=2, m_{jk}=1$, $q_{kk} \in \G_3$ or $q_{jj}=-1$, then
\begin{equation} \left[(\ad x_k)^2x_j, (\ad x_k)x_j \right]_c = 0; \label{Brelation}
\end{equation}
if there exist distinct $k, j, l$ satisfying $\ m_{kj}=2, m_{jk}=
m_{jl}=1$, $q_{kk} \in \G_3$ or $q_{jj}=-1$, then
\begin{equation}
\label{Brelation2} \left[ (\ad x_k)^2 (\ad x_j) x_l, (\ad x_k)x_j
\right]_c = 0;
\end{equation}
if $i,j$ determines a connected component of the Dynkin diagram of type $G_2$ and $q_{kk} \in \G_4$ or
$q_{jj}=-1$, then
\begin{align}
\left[ (\ad x_k)^3x_j, (\ad x_k)^2x_j \right]_c=&0, \label{G21'}
\\  \left[x_k,  \left[x_k^2x_jx_kx_j \right]_c \right]_c =&0, \label{G22'}
\\  \left[ \left[ x_k^2x_jx_kx_j \right]_c  , \left[x_kx_j\right]_c  \right]_c
=&0, \label{G23'}
\\ \left[  \left[x_k^{2}x_j\right]_c  , \left[ x_k^2x_jx_kx_j \right]_c \right]_c =&
0. \label{G24'}
\end{align}

Moreover, a basis of ${\mathfrak
B}(V)$ is given explicitly by:
$$x_{\beta_{1}}^{h_{1}} x_{\beta_{2}}^{h_{2}} \dots x_{\beta_{P}}^{h_{P}}, \qquad 0 \le h_{j} \le N_{\beta_j} - 1, \text{ if }\, \beta_j \in S_I, \quad  1\le
j \le P. \qed  $$
\end{theorem}

\subsection{Generalized Dynkin diagrams}

Given a braided vector space of diagonal type, with matrix $(q_{ij})_{1\leq i,j\leq\theta}$, there is a \textit{generalized Dynkin diagram} \cite{H2} associated to it, in such a way that two braided
vector spaces of diagonal type have the same generalized Dynkin diagram if and
only if they are twist equivalent. This diagram is a labeled graph
with vertices $1,\ldots,\theta$, each one labeled with the corresponding scalar $q_{ii}$. There is an edge between two different vertices $i$ and $j$ if $q_{ij}q_{ji}\neq1$ and it is labeled with this scalar.

The generalized Dynkin diagrams whose associated Nichols algebra is finite-dimensional were classified in \cite{H2}. We explicitly exhibit this classification for the case $A_2$. Let $q\in\G_N$, $N>2$, then the following are all the generalized Dynkin diagrams of standard type $A_2$:
\begin{align*}
\mD_1&=\xymatrix{
{\circ}_q \ar@{-}[r]^{q^{-1}} & {\circ}_{-1}}&&\mD_2=\xymatrix{
{\circ}_{-1} \ar@{-}[r]^{q} & {\circ}_{-1}}\\
\mD_3&=\xymatrix{
{\circ}_{-1} \ar@{-}[r]^{-1} & {\circ}_{-1}} && \mD_4=\xymatrix{
{\circ}_{q} \ar@{-}[r]^{q^{-1}} & {\circ}_{q}}
\end{align*}



\subsubsection{Yetter-Drinfeld modules} We denote by $\ydh$ the category of Yetter-Drinfeld
modules over $H$, where $H$ is a Hopf algebra with bijective
antipode. $M$ is an object of $\ydh$ if and only if there
exists an action $\cdot$ such that $(M,\cdot)$ is a (left)
$H$-module and a coaction $\delta$ such that $(M,\delta)$ is a
(left) $H$-comodule, subject to the following compatibility
condition:
\begin{equation*}
\delta(h\cdot m)=h_1m_{-1}S(h_3)\ot h_2\cdot m_{0}, \ \forall\, m\in
M, h\in H,
\end{equation*}
where $\delta(m)=m_{-1}\ot m_0$. Any $V\in \ydh$ becomes a braided vector space, \cite{Mo}.
If $\Gamma$ is a finite abelian group and
$H=\k \Gamma$, we denote $\ydg$ instead of $\ydh$. Any
$V\in\ydg$ is a braided vector space of diagonal type. Indeed, $V
= \oplus_{g\in \Gamma, \chi\in \VGamma}V_{g}^{\chi}$, where
$V_{g}^{\chi} = V^{\chi} \cap V_{g}$, $V_{g} = \{v\in V \mid
\delta(v) = g\otimes v\}$, $V^{\chi} = \{v\in V \mid  g \cdot v =
\chi(g)v \text{ for all } g \in \Gamma\}$.  The braiding is given
by $ c(x\otimes y) = \chi(g) y\otimes x$,   for all $x\in V_{g}$,
$g \in \Gamma$, $y\in V^{\chi}$, $\chi \in \VGamma$. Reciprocally, any braided vector space of diagonal type can be realized as a Yetter-Drinfeld module over the group algebra of an
abelian group.

\subsubsection{Nichols algebras}
If $V\in \ydh$, then the tensor algebra $T(V)$ admits a
unique structure of graded braided Hopf algebra in $\ydh$ such that
$V \subseteq \mP(V)$. The Nichols algebra $\bB(V)$ \cite{AS2} is the defined as the quotient of
$T(V)$ by the maximal element $I(V)$ of the class $\bS$ of all the homogeneous two-sided ideals $I \subseteq T(V)$
such that
\begin{itemize}
    \item $I$ is generated by homogeneous elements of degree $\geq
    2$,
    \item $I$ is a Yetter-Drinfeld submodule of $T(V)$,
    \item $I$ is a Hopf ideal: $\Delta(I) \subset I\ot T(V) +
    T(V)\ot I$.
\end{itemize}


\subsection{Lifting method}

Let $\Gamma$ be a finite group. The main steps of the \emph{Lifting Method} \cite{AS2} for
the classification of all finite-dimensional pointed Hopf algebras
with group of group-likes (isomorphic to) $\Gamma$ are:

\begin{itemize}
\item determine all $V\in\ydg$ such that the Nichols algebra $\B(V)$ is finite dimensional,

\item for such $V$, compute all Hopf algebras $H$ such
that $\gr H\simeq \B(V)\sharp \k\Gamma$. We call $H$ a
\emph{lifting} of $\B(V)$ over $\Gamma$.

\item Prove that any finite-dimensional pointed Hopf algebras with group $\Gamma$
is generated by group-likes and skew-primitives.
\end{itemize}

\section{Generation in degree one}\label{sec:gen in 1}

Throghout this Section, $\Gamma$ will denote a finite abelian group and $S= \bigoplus_{n \geq 0} S(n)$ a finite-dimensional graded braided Hopf algebra in $\gyd$ such that $S(0)=\k 1$. We fix a basis $\{x_1, \ldots, x_{\theta}\}$ of $V:=S(1)$, with $x_i \in S(1)^{\chi_i}_{g_i}$ for some $g_i \in \Gamma$ and $\chi_i \in \hat{\Gamma}$, and call $q_{ij}:=\chi_j(g_i)$.

We will show that given such $S$, if $V$ is a braided vector space of standard type, then $S$ is the Nichols algebra $\B(V)$ associated to $V$. In particular, we will obtain the main result of this work, that is that any finite-dimensional pointed Hopf algebra over $\Gamma$ with infinitesimal braiding of standard type is generated by group-like and skew-primitive elements.

First, we prove in the next Proposition that the quantum Serre relations hold in any such $S$, not necessarily of standard type. This result extends \cite[Lemma 5.4]{AS4}.

\begin{pro}\label{pro:qSerre} Let $S$ as above. Then,
\begin{equation}\label{qSerre}
    \ad_c(x_i)^{1+m_{ij}}(x_j)=0, \quad \mbox{for all }i \neq j \mbox{ such that }q_{ii}^{m_{ij}+1} \neq 1.
\end{equation}
\end{pro}
\pf
Suppose that $\ad_c(x_i)^{1+m_{ij}}(x_j) \neq 0$ for some $i \neq j$ such that $q_{ii}^{m_{ij}+1} \neq 1$ (so $q_{ii}^{m_{ij}}q_{ij}q_{ji}=1$ by definition of $m_{ij}$).

To start with, we begin as in \cite[Lemma 5.4]{AS4}. Set $m=m_{ij}$, $q=q_{ii}$, $y_1:=x_i$, $y_2:=x_j$ and $y_3:= \ad_c(x_i)^{1+m}(x_j)$. Also,
\begin{align*}
&h_1=g_i, &&h_2=g_j, &&h_3=g_i^{m+1}g_j, \\ &\eta_1=\chi_i, &&\eta_2=\chi_j, &&\eta_3= \chi_i^{m+1}\chi_j,
\end{align*}
so $y_k \in S^{\eta_k}_{h_k}, \ 1 \leq k \leq 3$. If $W=\k y_1+\k y_2+ \k y_3$, then $\bB(W)$ is finite-dimensional, because $y_3$ is a primitive element. Indeed, $W\subset \mP(S)$ hence we have a monomorphism $\B(W)\hookrightarrow S$. We compute the corresponding braiding matrix $\left(Q_{kl}= \eta_l(h_k) \right)_{1 \leq k,l \leq 3 }$, and consider the corresponding generalized Dynkin diagram:
\begin{equation}\label{diagram:Qserre}
\xymatrix{ & \circ^{q_{jj}} \ar@{-}[rd]^{q^{-m(m+1)}q_{jj}^{2}} &  \\  \circ^{q}\ar@{-}[ru]^{q^{-m}} \ar@{-}[rr]_{q^{m+2}} & & \circ^{q^{m+1}q_{jj}.} }
\end{equation}

In consequence, this diagram appears in \cite[Table 2]{H2}. We consider different cases.
\smallskip

\textbf{\emph{Case I:}} $Q_{kl}Q_{lk} \neq 1$ for all $1 \leq k<l \leq 3$.

By \cite[Lemma 9(ii)]{H2}, $1= \prod_{k<l} Q_{kl}Q_{lk}=q^{2-m(m+1)}q_{jj}^2$, and at least one of the vertices is labeled with $-1$. Notice that $q \neq -1$ because in such case $m=0$ (we assume $q^{m+1} \neq 1$). Also, $q_{jj} \neq q^{m+1}q_{jj}$ by hypothesis, so exactly one of the vertices is labeled with $-1$.
\begin{itemize}
  \item If $q_{jj}=-1$, then $1=(q^{m+1}q_{jj})(q^{-m(m+1)}q_{jj}^2)= -q^{1-m^2}$ and $m=1$ by the same Lemma, but this is a contradiction.
  \item If $q^{m+1}q_{jj}=-1$, then $1=qq^{m+2}=q^{m+3}$ and
  $$1=q_{jj}(q^{-m(m+1)}q_{jj}^{2})=q_{jj}^3q^{-m(m+3)+2m}=q_{jj}^3q^{2m}, $$
by the same Lemma, so
$$ -1= (-1)^3= q_{jj}^{3}q^{3m+3}= (q_{jj}^3q^{2m})q^{m+3}=q^{m+3}, $$
which is a contradiction. Therefore \eqref{diagram:Qserre} does not belong to this case.
\end{itemize}
\medskip

\textbf{\emph{Case II}:} $Q_{12}Q_{21}=q^{-m}=1$.

Here $m=0$, so we have
\begin{equation}\label{diagram:symmetriccase}
 \xymatrix{ \circ^{q} \ar@{-}[r]_{q^2} & \circ^{qq_{jj}} \ar@{-}[r]_{q_{jj}^2} & \circ^{q_{jj}}. }
\end{equation}
If $q_{jj}=-1$ then we have the connected subdiagram $ \xymatrix{ \circ^{q} \ar@{-}[r]_{q^2} & \circ^{-q} }$. Notice that this diagram has no vertices labeled with $-1$ and the labels of the vertices are different. Also the diagram is not of finite Cartan type and it does not correspond to the diagrams without $-1$ at the vertices in rows 5, 9, 11, 12, 15 of \cite[Table 1]{H2}, so we discard all of them.

If $q_{jj} \neq -1$ but $q=-1$ we have an analogous situation, so we consider also $q \neq -1$ and \eqref{diagram:symmetriccase} is a connected diagram of rank 3.

If $qq_{jj} \neq -1$, \cite[Lemma 9(i)]{H2} implies that one of the following holds:
\begin{itemize}
  \item it is of finite Cartan type, so it contains an $A_2$ Cartan subdiagram. Then $1= qq^2=(qq_{jj})q^2$ or $1= q_{jj}q_{jj}^2=(qq_{jj})q_{jj}^2$, so $q=1$ or $q_{jj}=1$;
  \item $q^3=1$, $q_{jj}, q_{jj}q \in \G_6 \cup \G_9$ and $q_{jj}q_{jj}^2=1$ or $q_{jj}^3=1$, $q, q_{jj}q \in \G_6 \cup \G_9$ and $qq^2=1$.
\end{itemize}
But neither of these cases is possible. In consequence, $qq_{jj} =-1$. Looking at \cite[Table 2]{H2} we see that $Q_{ii}Q_{i3}Q_{3i}=1$ for some $i \in \{1,2\}$ in all the cases. As both situations are analogous, we assume $i=1$: $q^3=1$. By \cite[Lemma 9(iii)]{H2}, one of the following is true:
\begin{itemize}
  \item $q_{jj}^3=1$, but $q_{jj}^3=-q^{-3}=-1$,
  \item $q_{jj}^4=1$,
  \item $q_{jj}=-q$.
\end{itemize}
We obtain a contradiction, so $m \neq 0$.
\medskip

\textbf{\emph{Case III:}} $Q_{13}Q_{31}=q^{m+2}=1$.

The corresponding diagram is: $$ \xymatrix{ \circ^{q} \ar@{-}[r]_{q^2} & \circ^{q_{jj}} \ar@{-}[r]_{q^{-2}q_{jj}^2} & \circ^{q^{-1}q_{jj}}. } $$
This diagram is analogous to \eqref{diagram:symmetriccase} exchanging $q_{jj}$ with $q_{jj}q^{-1}$ so we see that it does not belong to \cite[Table 2]{H2}. Therefore $q^{m+2} \neq 1$.
\medskip

\textbf{\emph{Case IV:}} $Q_{23}Q_{32}=1$. That is, $q_{jj}^2= q^{m(m+1)}$. We have the following diagram:
\begin{equation}\label{diagram:lastcase}
 \xymatrix{ \circ^{q_{jj}} \ar@{-}[r]_{q^{-m}} & \circ^{q} \ar@{-}[r]_{q^{m+2}} & \circ^{q^{m+1}q_{jj}}. }
\end{equation}
By the previous cases, this is a connected diagram of rank 3. As $m \neq 0$ and $q^{m+1} \neq 1$ we have $q \neq {-1}$. We analyze the different possibilities for the values on the vertices:

$\mathbf{q_{jj}=q^{m+1}q_{jj}=-1}$: In such case, $q^{m+1}=1$ and the diagram is $$ \xymatrix{ \circ^{-1} \ar@{-}[r]_{q} & \circ^{q} \ar@{-}[r]_{q} & \circ^{-1}, } $$ but it does not appear in Heckenberger's list.
\smallskip

$\mathbf{q_{jj}=-1, q^{m+1}q_{jj}\neq -1}$: By \cite[Table 2]{H2}, we have $1=Q_{22}Q_{23}Q_{32}=q^{m+3}$ and the diagram is
  $$ \xymatrix{ \circ^{-q^{-2}} \ar@{-}[r]_{q^{-1}} & \circ^{q} \ar@{-}[r]_{q^3} & \circ^{-1}. } $$
  Also, $ 1=q_{jj}^2= q^{m(m+1)}=q^{2m}=q^{-6} $. Notice that $q^3 \neq 1$ because $q^m \neq 1$, so $q \in \G_6$. But this diagram does not belong to Heckenberger's list.
\smallskip

$\mathbf{q_{jj}\neq -1, q^{m+1}q_{jj}=-1}$: As in the previous case, $1=Q_{22}Q_{21}Q_{12}=q^{1-m}$. By the definition of $m$ we conclude that $m=1$ and the diagram is the same as in the previous case, where again $q \in \G_6$ by the initial condition of case IV, and we have the same contradiction.
\smallskip

$\mathbf{q_{jj}, q^{m+1}q_{jj} \neq -1}$: By \cite[Lemma 9(i)]{H2}, one of the following holds:
\begin{itemize}
  \item it is of Cartan type. Therefore $q=q_{jj}$ and $m=1$, or $q=q^{m+1}q_{jj}=q^{-m-2}$. In both cases we arrive to the same diagram $$ \xymatrix{ \circ^{q} \ar@{-}[r]_{q^{-1}} & \circ^{q} \ar@{-}[r]_{q^3} & \circ^{q^3}. } $$
      It is easy to see that it is not of types $A_3, C_3$ because $q, q^2 \neq q^3$. But if it were of type $C_3$, $q=(q^3)^2=q^{-3}$, which is a contradiction.
  \item $q_{jj} \in \G_3$, $q\in \G_6 \cup \G_9$ and $1=q^{1-m}=q_{jj}q^{2m+3}$. Then $m=1$ and $q^5=q_{jj}^{-1}$, so $q^{15}=1$, but this is a contradiction with $q\in \G_6 \cup \G_9$.
  \item $q^{m+1}q_{jj} \in \G_3$, $q\in \G_6 \cup \G_9$ and $1=q_{jj}q^{-m}=q^{m+3}$. Again $q^{15}=1$, which is a contradiction with $q\in \G_6 \cup \G_9$.
\end{itemize}

In consequence \eqref{diagram:Qserre} is not of finite type, so we conclude the proof.
\epf

The following Lemmata show that if the braiding satisfies certain conditions related with braidings of standard type, then some extra relations hold in $S$. We consider the presentation of Nichols algebras of standard type from Theorem \ref{thm:presentacionNichols}. As this presentation is not minimal in some cases, we need first to discard some redundant relations in the next Lemma.

\begin{lem}\label{lem:redundantrelsB}
If there exist $j \neq k  \in \{1, \ldots, \theta \}$ such that $m_{kj}=1$, $m_{jk}=2$, but $q_{jj} \notin \G_3$ or $q_{kk} \neq -1$, then $\left[(\ad_c x_j)^2x_k, (\ad_c x_j)x_k \right]_c=0$.
\end{lem}

\pf
To prove this, by \cite[Lemma 5.5(ii)]{A} it is enough to consider two cases: $q_{jj} \in \G_3$, $q_{kk} \neq -1$, or $q_{jj} \notin \G_3$, $q_{kk}=-1$. In the first case we have
$$ x_j^3=0, \quad (\ad_c x_k)^2 x_j=x_k^2x_j-(1+q_{kk})q_{kj}x_kx_jx_k+q_{kk}q_{kj}^2x_jx_k^2=0. $$
In consequence we have $$ x_j^2x_kx_jx_k = (1+q_{kk})^{-1}q_{kj}^{-1}x_j^2x_k^2x_j, $$
and by \cite[Lemma 5.5(i)]{A} we conclude that $\left[(\ad_c x_j)^2x_k, (\ad_c x_j)x_k \right]_c=0$ (we can restrict to the Hopf subalgebra generated by $x_j,x_k$ in order to satisfy the conditions of such Lemma). The proof for the other case is analogous.
\epf

\begin{lem}\label{lemma:genqSerreA} Assume that there exist distinct $j, k, l \in \{1, \ldots, \theta \}$ such that $q_{kk}=-1$, $q_{kj}q_{jk}=q_{kl}^{-1}q_{lk}^{-1} \neq 1$, $q_{jl}q_{lj}=1$. Then,
\begin{equation}\label{genqSerreA}
    x_k^{2}=0, \quad \left[\ad_c x_j(\ad_c x_k(x_l)), x_k \right]_c=0.
\end{equation}
\end{lem}
\pf
The first equation follows easily because $x_k^{2}$ is primitive and the associated scalar is $1$. This implies that $(\ad_c x_k)^{2}x_j=(\ad_c x_k)^{2}x_l=0$.

For the second equation, we denote $u:= \left[\ad_c x_j(\ad_c x_k(x_l)), x_k \right]_c$, $g_u:=g_jg_k^{2}g_l \in \Gamma$, $\chi_u:=\chi_j\chi_k^{2}\chi_l\in \hat{\Gamma}$, $q:=q_{lk}q_{kl}$. By \cite[Lemma 5.8]{A}, $u$ is a primitive element.

We proceed as in the proof of the previous Lemma. Suppose that $u \neq 0$. Then the braiding of $y_1=x_j$, $y_2=x_k$, $y_3=x_l$ and $y_4=u$, with the corresponding elements $h_i \in \Gamma, \eta_i \in \hat{\Gamma}$, corresponds to one whose associated Nichols algebra is finite-dimensional. We obtain the following generalized Dynkin diagram attached to $(Q_{rs}=\eta_s(h_r))_{1 \leq r,s \leq 4}$:
\begin{equation}\label{diagram:genqSerreA}
    \xymatrix{ \circ^{q_{jj}} \ar@{-}[r]_{q^{-1}} \ar@{-}[d]_{q_{jj}^{2}q^{-2}} & \circ^{-1} \ar@{-}[d]_{q} \\ \circ^{q_{jj}q_{ll}} \ar@{-}[r]_{q_{ll}^{2}q^2} & \circ^{q_{ll}}.  }
\end{equation}
Notice that $q=-1$ implies that \eqref{diagram:genqSerreA} contains \eqref{diagram:symmetriccase} as a subdiagram, which is a contradiction to the finite dimension of the associated Nichols algebra. Therefore $q \neq -1$ and as any such diagram contains a 4-cycle, by \cite[Lemma 12]{H2} we have $q_{jj}^2q^{-2}=1$ or $q_{ll}^2q^2=1$. By the symmetry of the diagram we can assume $q_{jj}= \pm q$.

If also $q_{ll}= \pm q^{-1}$, since $Q_{44}=q_{jj}q_{ll} \neq 1$, then the diagram is of the form $\xymatrix{\circ^{q} \ar@{-}[r]_{q^{-1}} & \circ^{-1} \ar@{-}[r]_{q} & \circ^{-q^{-1}} }$. But this is a contradiction with \cite[Lemma 9(iii)]{H2}. Therefore we have a connected diagram of rank 4:
$$ \xymatrix{\circ^{\pm q} \ar@{-}[r]_{q^{-1}} & \circ^{-1} \ar@{-}[r]_{q} & \circ^{q_{ll}}  \ar@{-}[r]_{q_{ll}^2q^2} & \circ^{q_{jj}q_{ll}}. } $$
Suppose that $q_{jj}=-q$. As $Q_{11}Q_{12}Q_{21} \neq 1$, we deduce from \cite[Table 3]{H2} that $$0=(1-Q_{11}^3)(Q_{11}^2Q_{12}Q_{21}-1)=(1+q^3)(q-1),$$
but we discard this case by \cite[Table 3]{H2}.

Therefore $q_{jj}=q$. We obtain that there are no diagrams in \cite[Table 3]{H2} such that $Q_{22}=-1$, $Q_{11}=Q_{44}Q_{33}^{-1}=q \neq \pm1$, so the above diagram does not belong to such list. Therefore $u=0$.
\epf

\begin{lem}\label{lemma:genqSerreB} Assume that there exist $j \neq k  \in \{1, \ldots, \theta \}$ such that $m_{kj}=1$, $m_{jk}=2$.
\smallskip

\emph{(a)} If $q_{jj} \in \G_3$ and $q_{kk}=-1$, then the following relation holds:
\begin{equation}\label{genqSerreB2}
    \left[(\ad_c x_j)^2x_k, (\ad_c x_j)x_k \right]_c=0.
\end{equation}

\emph{(b)} If $V$ is standard and there exist $l \neq j,k  \in \{1, \ldots, \theta \}$ such that $m_{jl}=m_{lj}=0$, $m_{kl}=1$ and $(1+q_{kk})(1-q_{jj}^3)=0$, then:
\begin{equation}\label{genqSerreB3}
    \left[(\ad_c x_k)^2(\ad_c x_j)x_l, (\ad_c x_j)x_k \right]_c=0.
\end{equation}

\end{lem}

\pf
(a) We proceed as in the previous Lemmata. Assume that $v:=\left[(\ad_c x_j)^2x_k, (\ad_c x_j)x_k \right]_c \neq 0$. By \cite[Lemma 5.9]{A}, $v$ is a primitive element: notice that $x_j^3=0$, or $x_k^2=0$, or $q_{jj}^2q_{jk}q_{kj}=q_{kk}q_{jk}q_{kj}=1$  because $S$ is finite dimensional.

Call $y_1=x_j$, $y_2=x_k$, $y_3=v$, and $h_i \in \Gamma$, $\eta_i \in \hat{\Gamma}$, $i=1,2,3$ the corresponding elements. In consequence, the braiding matrix $(Q_{rs}=\eta_s(h_r))_{1 \leq r,s \leq 3}$ corresponds to one in Heckenberger's list. The associated generalized Dynkin diagram is
$$ \xymatrix{ \circ^{q_{jj}} \ar@{-}[rr]_{q} \ar@{-}[rd]_{q^2} &  & \circ^{-1} \\ & \circ^{q^6} \ar@{-}[ru]_{q^3} & }, \quad q:=q_{jk}q_{kj}. $$
As the diagram is finite, $Q_{33}=q^{6} \neq 1$, but then this contradicts \cite[Lemma 9(iii)]{H2}. Therefore, $v=0$.
\medskip

(b) By the previous item, Lemma \ref{lem:redundantrelsB} and \cite[Lemma 5.9(b)]{A}, $$w:= \left[(\ad_c x_k)^2(\ad_c x_j)x_l, (\ad_c x_j)x_k \right]_c$$ is a primitive element. If we suppose that $w\neq 0$, we work as in previous cases for each possible diagram calling $y_1=x_j$, $y_2=x_k$, $y_3=x_l$, $y_4=w$, and $h_i \in \Gamma$, $\eta_i \in \hat{\Gamma}$, $i=1,2,3,4$ the corresponding elements: the braiding matrix $(Q_{rs}=\eta_s(h_r))_{1 \leq r,s \leq 3}$ corresponds to one in Heckenberger's list.

\begin{itemize}
  \item $q_{kk}=-1, \ q_{jj}^2q_{kj}q_{jk}=1=q_{kj}q_{jk}q_{kl}q_{lk}$: the corresponding diagram for $(Q_{rs})$ is
  $$ \xymatrix{ \circ^{q} \ar@{-}[r]^{q^{-2}} \ar@{-}[rd]_{q^2} & \circ^{-1} \ar@{-}[r]^{q^2}\ar@{-}[d]^{q^{-4}} & \circ^{q_{ll}} \ar@{-}[ld]^{q^4q_{ll}^2}  \\ & \circ^{qq_{ll}}  &  }, \qquad q:=q_{jj}. $$
  By \cite[Lemma 9(ii)]{H2} $q^4=1$. Then the vertices 1,3,4 determine a diagram of type \eqref{diagram:symmetriccase}, which is not in Heckenberger's list.
  \item $\ q_{jj}^2q_{kj}q_{jk}=q_{kj}q_{jk}=q_{kk}q_{kl}q_{lk}=1$: the diagram for this case is
  $$ \xymatrix{ \circ^{q} \ar@{-}[r]^{q^{-2}} \ar@{-}[rd]_{q^2} & \circ^{q^2} \ar@{-}[r]^{q^{-2}} & \circ^{q_{ll}} \ar@{-}[ld]^{q^{-4}q_{ll}^2}  \\ & \circ^{qq_{ll}}  &  }, \qquad q:=q_{jj}\in \G_3 \cup \G_4. $$
  If $q^4=1$, again we have \eqref{diagram:symmetriccase} as subdiagram. If $q \in \G_3$ we have $q_{ll}= \pm q^2$, because there are no 4-cycles. As $Q_{44} \neq 1$, we should have $q_{ll}=-q^2$, but in such case we have a connected diagram of rank 4 with $m_{32}=4$, because $Q_{33}=-q^{2}$ and $Q_{23}Q_{32}=q^2 \in \G_3$, which is a contradiction.
  \item $q_{kk}=-1, \ q_{jj}=-q_{kj}q_{jk}\in G_3, \ q_{kj}q_{jk}q_{kl}q_{lk}=1$: the diagram is
  $$ \xymatrix{ \circ^{q} \ar@{-}[r]^{-q} \ar@{-}[rd]_{q^2} & \circ^{-1} \ar@{-}[r]^{-q^2}\ar@{-}[d]^{q^2} & \circ^{q_{ll}} \ar@{-}[ld]^{qq_{ll}^2}  \\ & \circ^{qq_{ll}}  &  }, \qquad q:=q_{jj} \in \G_3. $$
  This is not of finite type by \cite[Lemma 9(ii)]{H2}.
  \item $q_{jj}=-q_{kj}q_{jk} \in \G_3, q_{kk}q_{kj}q_{jk}=q_{kk}q_{kl}q_{lk}=1$: now the diagram is
$$ \xymatrix{ \circ^{q} \ar@{-}[r]^{-q} \ar@{-}[rd]_{q^2} & \circ^{-q^2} \ar@{-}[r]^{-q} & \circ^{q_{ll}} \ar@{-}[ld]^{q^2q_{ll}^2}  \\ & \circ^{qq_{ll}}  &  }, \qquad q:=q_{jj} \in \G_3. $$
First $q_{ll}= \pm q^2$ because there are no 4-cycles, and second $Q_{44} \neq 1$ so $q_{ll}=-q^2$. Transforming the diagram by the symmetry at vertex 4, it is Weyl equivalent to
$$ \xymatrix{ \circ^{-1} \ar@{-}[r]^{q} & \circ^{1} \ar@{-}[r]^{-q} & \circ^{-q^2} \ar@{-}[r]^{-q} & \circ^{q_{ll}}, } $$
whose associated Nichols algebra is not finite-dimensional, a contradiction.
\end{itemize}
In all the cases we obtain a contradiction, so $w=0$.
\epf

We are now able to prove the main results of this Section: Theorems \ref{thm:generationdegree} and \ref{thm:generationdegree1}.

\begin{theorem}\label{thm:generationdegree} Let $S= \oplus_{n \geq 0} S(n)$ be a finite-dimensional graded Hopf algebra in $\gyd$, $\Gamma$ a finite abelian group, such that $S(0)=\k 1$. Fix a basis $x_1, \ldots, x_{\theta}$ of $V:=S(1)$, with $x_i \in S(1)^{\chi_i}_{g_i}$ for some $g_i \in \Gamma$ and $\chi_i \in \hat{\Gamma}$, and call $q_{ij}:=\chi_j(g_i)$. Assume that
\begin{itemize}
  \item S is generated as an algebra by $S(0)\oplus S(1)$, and
  \item $V$ is a standard braided vector space.
\end{itemize}
Then $S \cong \bB(V)$.
\end{theorem}
\pf
The canonical surjection $T(V) \twoheadrightarrow \bB(V)=T(V)/I(V)$ induces a surjection
$$ \pi: S \twoheadrightarrow \bB(V), $$
of braided graded Hopf algebras, so we can consider $S=T(V)/I$, for some graded braided Hopf ideal $I$ of $T(V)$, generated in degree $\geq 2$, $I \subseteq I(V)$.

Suppose that $I(V) \supsetneqq I$. Then at least one generator of $I(V)$ as in Theorem \ref{thm:presentacionNichols} is not in $I$. Consider a generator $\mathbf{x} \in I(V) \setminus I$ of minimal degree $k$. Then
$$ \Delta(\mathbf{x})= \mathbf{x} \otimes 1+1 \otimes \mathbf{x}+ \sum_{j=1}^n b_j \otimes c_j \in I(V) \otimes T(V) + T(V) \otimes I(V), $$
for some homogeneous elements $b_j, c_j \in \bigoplus_{i=1}^{k-1} T^i(V)$, satisfying $\deg(b_j)+\deg(c_j) =k$. We can consider for each $j$ that $b_j \in I(V)$ or $c_j\in I(V)$. If $b_j\in I(V)$, then it is a linear combination of elements $a\mathbf{y}b$, with $a,b \in T(V)$ and $\mathbf{y}$ a generator of I(V), and as $$ \deg(\mathbf{y}) = \deg (\mathbf{x})-\deg(a)-\deg(b)-\deg(c_j) < \deg(\mathbf{x}) =k,$$
we have $\mathbf{y}\in I$. Therefore $b_j \in I$. The same holds if $c_j \in I(V)$. In consequence, $\mathbf{x}$ is primitive in $S$, because $\pi$ is a morphism of braided Hopf algebras.

By Proposition \ref{pro:qSerre} and Lemmata \ref{lem:redundantrelsB}, \ref{lemma:genqSerreA} and \ref{lemma:genqSerreB}, we have $\mathbf{x}=x_\alpha^{N_{\alpha}}$ for some $\alpha \in \Delta_+$, or there exist $j \neq k \in \{1, \ldots , \theta \}$ such that $m_{jk}=3$, $m_{kj}=1$, $(1-q_{jj}^4)(1+q_{kk})=0$, $\left( \begin{array}{cc} q_{jj} & q_{jk} \\ q_{kj} & q_{kk} \end{array} \right)$ is a standard braiding of type $G_2$ and $\mathbf{x}=[u]_c$, see \ref{subsec:braideddiagonal}, where
$$ u \in \{ x_j^3x_kx_jx_k, \ \ x_j^3x_kx_j^2x_k,  \ \  x_j^2x_kx_jx_kx_jx_k,  \ \  x_j^2x_kx_j^2x_kx_jx_k \}. $$
Call $g_{\mathbf{x}} \in \Gamma$, $\chi_{\mathbf{x}} \in \hat{\Gamma}$ the associated elements. We discard easily the case $\mathbf{x}=x_\alpha^{N_{\alpha}}$, because in such case $$\chi_{\mathbf{x}}(g_{\mathbf{x}})=q_{\alpha}^{N_{\alpha}^2}=1,$$
($\ord(q_{\alpha})=N_{\alpha}$) and $S$ is finite-dimensional.

Suppose $\mathbf{x}=[u]_c$. Call $\eta_1=\chi_j$, $\eta_2=\chi_k$, $\eta_3=\chi_{\mathbf{x}}$, $h_1=g_j$, $h_2=g_k$, $h_3=g_{\mathbf{x}}$. As in the proof of previous Lemmata, the braiding corresponding to the matrix $(Q_{rs}=\eta_s(h_r))_{1 \leq r,s \leq 3}$ appears in Heckenberger's list. The possible diagrams for the vertices $j,k$ are
\begin{itemize}
  \item $\xymatrix{ \circ^{\zeta} \ar@{-}[r]^{\zeta} & \circ^{\zeta^3} }$, $\zeta \in \G_4$;
  \item $\xymatrix{ \circ^{\zeta} \ar@{-}[r]^{-1} & \circ^{-1} }$, $\zeta \in \G_6$;
  \item $\xymatrix{ \circ^{\zeta^2} \ar@{-}[r]^{\zeta} & \circ^{\zeta^7} }$, $\zeta \in \G_8$;
  \item $\xymatrix{ \circ^{\zeta^2} \ar@{-}[r]^{\zeta^3} & \circ^{-1} }$, $\zeta \in \G_8$;
  \item $\xymatrix{ \circ^{\zeta} \ar@{-}[r]^{\zeta^5} & \circ^{-1} }$, $\zeta \in \G_8$.
\end{itemize}
With three exceptions we conclude all of possible pairs of braidings and $u$ give diagrams not in Heckenberger's list because $Q_{33}=1$ or
\begin{align*}
 Q_{12}Q_{21}\neq1, && Q_{13}Q_{31}\neq1 && \text{and} && Q_{23}Q_{32} \neq 1
\end{align*}
and thus it is a triangle but $\prod_{1 \leq r<s \leq 3}Q_{rs}Q_{sr}\neq1.$

The remaining cases are:
\begin{enumerate}
\item $\xymatrix{ \circ^{\zeta} \ar@{-}[r]^{\zeta} & \circ^{\zeta^3} }$, \qquad $u=x_j^3x_kx_j^2x_k$,
\item $\xymatrix{ \circ^{\zeta} \ar@{-}[r]^{-1} & \circ^{-1} }$,\qquad $u=x_j^3x_kx_jx_k$,
\item $\xymatrix{ \circ^{\zeta} \ar@{-}[r]^{-1} & \circ^{-1} }$, \qquad $u=x_j^2x_kx_jx_kx_jx_k$,
\end{enumerate}
and corresponding diagrams of $(Q_{rs})$ are:
\begin{enumerate}
\item $\xymatrix{ \circ^{\zeta} \ar@{-}[r]^{\zeta} & \circ^{\zeta^3} \ar@{-}[r]^{\zeta} & \circ{\zeta^3} }$,
\item $\xymatrix{ \circ^{\zeta^4} \ar@{-}[r]^{\zeta^2} & \circ^{\zeta} \ar@{-}[r]^{-1} & \circ^{-1}}$,
\item $\xymatrix{ \circ^{\zeta} \ar@{-}[r]^{\zeta^5} & \circ^{\zeta} \ar@{-}[r]^{-1} & \circ^{-1} }$,
\end{enumerate}
but they are Cartan braidings associated to non-finite Cartan matrices, which is a contradiction.
\epf

The following Theorem is in agreement with Conjecture \cite[Conj. 1.4]{AS2}. As braidings of standard type properly include those of finite Cartan type, this result extends \cite[Thm. 5.5]{AS4}.

\begin{theorem}\label{thm:generationdegree1}
Let $H$ be a finite-dimensional pointed Hopf algebra over an abelian group $\Gamma$ such that its infinitesimal braiding is of standard type. Then $H$ is generated by its group-like and skew-primitive elements.
\end{theorem}

\pf
Let $\gr H=R\#\k \Gamma$, $V=R(1)$. Then $H$ is generated by its group-like and skew-primitive elements if and only if $R$ is the Nichols algebra $\B(V)$. Let $S$ be the graded dual $R^*$ in $\ydg$. Notice that $S(1)=R(1)^*$ has the same braiding as $R(1)$. By \cite[Lem. 5.5]{AS2} it is enough to show that $S$ is a Nichols algebra. This follows by Theorem \ref{thm:generationdegree}.
\epf

\section{Liftings of the quantum Serre relations}

Let $\B$ a finite dimensional Nichols algebra of standard type, with braiding $(q_{ij})_{1\leq i,j\leq \theta}$, $\theta$ the rank of $\B$. Let $H$ be a pointed Hopf algebra over an abelian group $\Gamma$ such that $\gr H=\B\#\k\Gamma$. In this Section we show that quantum Serre relations in $\B$ are lifted in $H$ as elements in $\k\Gamma$. We include this result in this work since its proof heavily resembles the one of Proposition \ref{pro:qSerre}. Moreover, we distinguish those cases in which these relations can only be lifted as zero.

Let $m_{ij}$ be as in \ref{subsec:bvs}. For $1\leq i\neq j\leq \theta$, set
$$
\chi_{ij}=\xi_i^{m_{ij}+1}\chi_j, \quad g_{ij}=g_i^{m_{ij}+1}g_j.
$$

\begin{lem}\label{lem:lifting quantum}
Let $1\leq i\neq j\leq \theta$. Assume $q_{ii}^{m_{ij}+1}\neq 1$. Then
\begin{equation}\label{eqn:chilgl}
(\chi_{ij},g_{ij})\neq (\chi_l,g_l), \quad \forall\,1\leq l\leq \theta.
\end{equation}
\end{lem}
\begin{proof}
Assume there exists $l$ such that \eqref{eqn:chilgl} holds. Then it follows as in the proof of Proposition \ref{pro:qSerre} that either $l=i$ or $l=j$, since otherwise we would get a subdiagram of the diagram $\mD$ related to the braiding which would not be of finite type.

Now, $l\neq j$, since otherwise $\chi_{ij}(g_{ij})=q_{ii}^{m_{ij}+1}q_{jj}=\chi_j(g_j)=q_{jj}$ and $q_{ii}^{m_{ij}+1}\neq 1$. But if $l=i$, we would have $$q_{ii}^2=\chi_i(g_i)^2=\chi_i(g_{ij})\chi_{ij}(g_i)=q_{ii}^{2(m_{ij}+1)}q_{ji}q_{ij}=q_{ii}^{m_{ij}+2},$$ since, as $(q_{ij})$ is of standard type and $q_{ii}^{m_{ij}+1}\neq 1$, $q_{ii}^{m_{ij}}q_{ij}q_{ji}=1$. Therefore, by definition of $m_{ij}$, we have $m_{ij}=0$. In this case, $q_{ii}=\chi_i(g_i)=\chi_{ij}(g_{ij})=q_{ii}q_{jj}$, a contradiction.
\end{proof}

Let $a_i\in\mP_{1,g_i}^{\chi_i}(H)$ such that $a_i$ is mapped to $x_i\in\B$ via $A_1\twoheadrightarrow A_1/A_0\cong\B(1)$. In particular, it follows that $\ad(a_i)^{m_{ij}+1}(a_j)\in\mP_{g_{ij}}^{\chi_{ij}}(H)$.

If $1\leq i\neq j\leq \theta$, denote by $\mD_{ij}$ the subdiagram of $\mD$ with vertices $i,j$ and by $Q_{ij}$ the corresponding submatrix of $(q_{kl})_{1\leq k,l\leq \theta}$. 

\begin{pro}\label{pro:lifting quantum}
Let $H,a_i$ as above, $1\leq i\neq j\leq \theta$. Assume $q_{ii}^{m_{ij}+1}\neq 1$. Then there exists $\lambda\in\k$ such that
$$
\ad(a_i)^{m_{ij}+1}(a_j)=\lambda(1-g_i^{m_{ij}+1}g_j).
$$
Moreover, $\lambda$ can be non-zero only in the following cases:
\begin{enumerate}
\item $m_{ij}=3$ and
\begin{itemize}
\item[(i)] $\mD_{ij}=\xymatrix{
{\circ}_q \ar@{-}[r]^{q^{-3}} & {\circ}_{q^3}}$, $q\in\mathbb{G}_7$ and $Q_{ij}=\left(\begin{smallmatrix}
                                                                                q&q^3\\q&q^3
                                                                                \end{smallmatrix}\right),
$
\item[(ii)] $\mD_{ij}=\xymatrix{
{\circ}_\xi \ar@{-}[r]^{\xi^{-3}} & {\circ}_{-1}}$, $\xi\in\mathbb{G}_8$ and $Q_{ij}=\left(\begin{smallmatrix}
                                                                                \xi&-1\\ \xi&-1
                                                                                \end{smallmatrix}\right).
$
\end{itemize}
\item $m_{ij}=2$ and
\begin{itemize}
\item[(i)] $\mD_{ij}=\xymatrix{
{\circ}_q \ar@{-}[r]^{q^{-2}} & {\circ}_{q^2}}$, $q\in\mathbb{G}_5$ and $Q_{ij}=\left(\begin{smallmatrix}
                                                                                q&q^2\\q&q^2
                                                                                \end{smallmatrix}\right),
$
\item[(ii)] $\mD_{ij}=\xymatrix{
{\circ}_q \ar@{-}[r]^{q^{-2}} & {\circ}_{-1}}$, $q\in\mathbb{G}_6$ and $Q_{ij}=\left(\begin{smallmatrix}
                                                                                q&-1\\q&-1
                                                                                \end{smallmatrix}\right).
$
\end{itemize}
\item $m_{ij}=1$ implies
\begin{itemize}
\item[(i)] $\mD_{ij}=\xymatrix{
{\circ}_{q^m} \ar@{-}[r]^{q^{-m}} & {\circ}_{q}}$, $q\in\mathbb{G}_{2m+1}$ and $Q_{ij}=\left(\begin{smallmatrix}
                                                                                q^m&q\\q^m&q
                                                                                \end{smallmatrix}\right),
$
\item[(ii)] $\mD_{ij}=\xymatrix{
{\circ}_{q} \ar@{-}[r]^{-q} & {\circ}_{-1}}$, $q\in\mathbb{G}_4$ and $Q_{ij}=\left(\begin{smallmatrix}
                                                                                q&-1\\q&-1
                                                                                \end{smallmatrix}\right),$
\item[(iii)] $\mD_{ij}=\xymatrix{
{\circ}_{-\xi} \ar@{-}[r]^{-\xi^{-1}} & {\circ}_{\xi}}$ $\xi\in\mathbb{G}_3$, and $Q_{ij}=\left(\begin{smallmatrix}
                                                                                -\xi&\xi\\-\xi&\xi
                                                                                \end{smallmatrix}\right),$
\item[(iv)] $\mD_{ij}=\xymatrix{
{\circ}_{q} \ar@{-}[r]^{q^{-1}} & {\circ}_{q^2}}$, $q\in\mathbb{G}_8$ and $Q_{ij}=\left(\begin{smallmatrix}
                                                                                q&q^2\\ q&q^2
                                                                                \end{smallmatrix}\right).
$
\end{itemize}
\item $m_{ij}=0$ implies
\begin{itemize}
\item[(i)] $\mD_{ij}=\xymatrix{
{\circ}_{q} \ar@{}[r] & {\circ}_{q^{-1}}}$, $q\in \G_N$, $N>1$ and $Q_{ij}=\left(\begin{smallmatrix}
                                                                                q&q^{-1}\\ q&q^{-1}
                                                                                \end{smallmatrix}\right).
$
\end{itemize}
\end{enumerate}
\end{pro}
\begin{proof}
By \cite[Lemma 5.4]{AS1} we know that $\mP_{g,1}^\eps(H)=\k(1-g)$ and that if $\chi\neq \eps$ then $\mP_{g,1}^\chi(H)\neq 0$ if and only if there is $1\leq l\leq \theta$ such that $(g_{ij},\chi_{ij})=(g_l,\chi_l)$. As $\ad(a_i)^{m_{ij}+1}(a_j)\in\mP_{g_{ij}}^{\chi_{ij}}(H)$ then the first part of the Proposition follows from Lemma \ref{lem:lifting quantum}. For the second part, as $\lambda$ can be chosen $\neq0$ only when  $\chi_i^{m_{ij}+1}\chi_j=\eps$, the Proposition follows by evaluating $\chi_{ij}$ in $g_i$ and $g_j$ to determine when $q_{ii}^{m_{ij}+1}q_{ij}=q_{ji}^{m_{ij}+1}q_{jj}=1$, taking into account that $q_{ii}^{m_{ij}}q_{ij}q_{ji}=1$. We fully develop the case $m_{ij}=3$ as an example. As $q_{ii}^{m_{ij}+1}\neq 1$ must hold, two cases are left, namely those corresponding to the diagrams in (1)(i) and (ii) of the Proposition. Let $(q_{ij})_{1\leq i,j\leq\theta}$ be the braiding. In the first case, we have $q_{ii}=q, q_{jj}=q^3$ and $q_{ij}q_{ji}=q^{-3}$. Then
$$
\chi_i^4\chi_j(g_i)=q^4q_{ij}, \quad \chi_i^4\chi_j(g_j)=q_{ji}^4q^3.
$$
Then, if $\chi_{ij}=\eps$, we have $1=q^4q_{ij}q_{ji}^4q^3=q^4q_{ji}^3$ and $1=q_{ji}^4q^3=q_{ji}q^{-1}$. Therefore, $q_{ji}=q$, $q^7=1$, $q_{ij}=q^{-4}=q^3$.

In the second case, we have $q_{ii}=\xi, q_{jj}=-1$, $q_{ij}q_{ji}=\xi^{-3}$, $\xi\in\mathbb{G}_8$. Then, if $\xi_{ij}=\eps$, $1=\xi^4q_{ij}=-q_{ji}^4$. Then $q_{ij}=-1$ and $q_{ji}=-\xi^{-3}=\xi$.
\end{proof}
%
%
%
%
%
%
%
%


\end{document}